\documentclass[11pt]{article}

\usepackage{amsmath,amssymb,amsfonts,amsbsy,mathrsfs}

\renewcommand{\subsection}{\subsubsection}

\begin{document}

\title{\bf Shock waves and characteristic discontinuities in ideal compressible two-fluid MHD}

\author{{\bf Lizhi Ruan}\\
Hubei Key Laboratory of Mathematical Physics, \\ School of Mathematics and Statistics, \\
Central China Normal University, Wuhan 430079, P.R. China \\
E-mail: rlz@mail.ccnu.edu.cn
\and
{\bf Yuri Trakhinin}\\
Sobolev Institute of Mathematics, Koptyug av. 4, 630090 Novosibirsk, Russia\\
and\\
Novosibirsk State University, Pirogova str. 2, 630090 Novosibirsk, Russia\\
E-mail: trakhin@math.nsc.ru
}

\date{ }

\maketitle

\begin{abstract}
\noindent We are concerned with a model of ideal compressible isentropic two-fluid magnetohydrodyna\-mics (MHD). Introducing an entropy-like function, we reduce the equations of two-fluid MHD to a symmetric form which looks like the classical MHD system written in the nonconservative form in terms of the pressure, the velocity, the magnetic field and the entropy. This gives a number of instant results. In particular, we conclude that all compressive extreme shock waves exist locally in time in the limit of weak magnetic field. We write down a condition sufficient for the local-in-time existence of current-vortex sheets in two-fluid flows. For the 2D case and a particular equation of state, we make the conclusion that contact discontinuities in two-fluid MHD flows exist locally in time provided that the Rayleigh-Taylor sign condition on the jump of the normal derivative of the pressure is satisfied at the first moment.
\end{abstract}



\section{Introduction}
\label{s1}

We consider the following equations of ideal compressible inviscid isentropic two-fluid magnetohydrodynamics (MHD):
\begin{equation}\label{1}
\left\{
\begin{array}{l}
 \partial_tn  +{\rm div}\, (n u )=0,\\[6pt]
  \partial_t\rho  +{\rm div}\, (\rho u )=0,\\[6pt]
 \partial_t((\rho +n) u ) +{\rm div}\,((\rho +n)u\otimes u -{H}\otimes{H} ) +
{\nabla}q=0, \\[6pt]
 \partial_t{H} -{\nabla}\times (u {\times}{H})=0,
\end{array}
\right.
\end{equation}
where $u\in\mathbb{R}^3$ denotes the two-fluid flow velocity, $n$ and $\rho$ are the densities of two fluids, $H\in\mathbb{R}^3$ is the magnetic field, $q =P+\frac{1}{2}|{H} |^2$ is the total pressure, and
\begin{equation}
P=P(\rho ,n)= \rho^{\alpha} +An^{\gamma}
\label{2}
\end{equation}
denotes the pressure for some constants $\alpha \geq 1$, $\gamma \geq 1$ and $A>0$. System \eqref{1} is supplemented by the divergence constraint
\begin{equation}
{\rm div}\, {H} =0
\label{3}
\end{equation}
on initial data.

System \eqref{1} with $\alpha =1$ and $\gamma >1$ is the inviscid version of the viscous equations from \cite{Wen18} which are formally derived from the Vlasov-Fokker-Planck/compressible magnetohydrodynamics equations by applying ideas of Carrillo and Goudon in \cite{Carrillo06}. The coupled system of Vlasov-Fokker-Planck and compressible magnetohydrodynamics equations describes the motions of uncharged particles in a viscous inhomogeneous compressible conducting fluid whose global well-posedness and the large time behavior of classical solutions are proved in \cite{Jiang17}. That is, for $\alpha =1$ and $\gamma >1$, we can treat \eqref{1} as the system modelling the motion of the mixture of an inviscid fluid and particles in magnetic field, where $n$ is the density of the fluid and $\rho$ is the density of particles in the mixture. In this case, it is better to call \eqref{1} a two-phase MHD model.

The general case of the equation of state \eqref{2} when $\alpha \geq 1$ and  $\gamma \geq 1$ is considered, for example in \cite{Vasseur17}, for the viscous compressible two-fluid model without magnetic field which was, in particular, derived from physical considerations in \cite{Ishii75} (see also \cite{Vasseur17} for other references). System \eqref{1} with the equation of state \eqref{2} is the inviscid version of this model completed (as in \cite{Wen18}) by the influence of magnetic field.


Stability/existence results for the multidimensional case (2D and/or 3D) for shock waves and characteristic discontinuities (current-vortex sheets, contact and Alfv\'{e}n discontinuities) in ideal compressible one-fluid (classical) MHD  can be found in \cite{BThand,ChW1,ChW2,IT,MTT1,MTT2,T,T05,T09,Wang13}. For vortex sheets for a compressible inviscid liquid-gas two-phase model (without magnetic field) we can refer to the stability/existence results in \cite{HWY,RWWZ}. An essential improvement of the last results for the liquid-gas two-phase model was done in our very recent work \cite{RT1}, where vortex sheets were also considered for model \eqref{1} without magnetic field. In \cite{RT1} we also prove the local-in-time existence of all compressive shock waves in the mentioned liquid-gas model as well as, for example, in model \eqref{1} without magnetic field.

In this paper, using the crucial simple idea from our recent study in \cite{RT1} and introducing an entropy-like function, we reduce the equations \eqref{1} of two-fluid MHD to a symmetric form which looks like the classical MHD system written in the nonconservative form in terms of the pressure, the velocity, the magnetic field and the entropy. This gives a number of instant results. In particular, we conclude that all compressive extreme shock waves exist locally in time in the limit of weak magnetic field. We write down a condition sufficient for the local-in-time existence of current-vortex sheets in two-fluid flows. For the 2D case and when $\alpha =\gamma$ in \eqref{2}, we make the conclusion that contact discontinuities in two-fluid MHD flows exist locally in time provided that the Rayleigh-Taylor sign condition on the jump of the normal derivative of the pressure is satisfied at the first moment.

The rest of the paper is organized as follows. In Section 2, we derive the mentioned symmetric form of equations \eqref{1} as well as write down a so-called secondary symmetrization which was first proposed in \cite{T05} for classical MHD. In Section 3, we introduce all four types of strong discontinuities (shock waves, current-vortex sheets, contact and Alfv\'{e}n discontinuities) in two-fluid MHD. In Section 4, we discuss structural stability conditions for the initial data for the free boundary problems for all the types of strong discontinuities in two-pase MHD introduced in Section \ref{s3}. These conditions should guarantee the local-in-time existence in Sobolev spaces of solutions of these problems (usually they also provide the uniqueness of a solution).

\section{Symmetrizations of the two-fluid MHD equations}
\label{s2}

Let
\begin{equation}
n >0\quad\mbox{and}\quad\rho \geq 0.
\label{4}
\end{equation}
Using the crucial simple idea from \cite{RT1}, we introduce the entropy-like function
\[
S=\frac{\rho}{n}\geq 0.
\]
It indeed plays the role of entropy and we will call it the ``entropy'' because, as for the usual entropy, we get the equation
\[
\frac{{\rm d} S}{{\rm d}t}=0,
\]
with ${\rm d} /{\rm d} t =\partial_t+({u} \cdot{\nabla} )$, following from the first two equations of system \eqref{1}. Moreover, introducing the total density $R=\rho +n >0$, we can equivalently rewrite the first two equations of  \eqref{1} as
\begin{equation}
\frac{1}{R}\,\frac{{\rm d} R}{{\rm d}t} +{\rm div}\,u =0,\quad \frac{{\rm d} S}{{\rm d}t}=0.
\label{5}
\end{equation}

We now recalculate the pressure $P$ in terms of $R$ and $S$. We have
\[
n=\frac{R}{S+1},\quad \rho =\frac{RS}{S+1}.
\]
Then
\[
P=P(R,S) =\left( \frac{RS}{S+1}\right)^{\alpha} + A\left( \frac{R}{S+1}\right)^{\gamma} >0
\]
and, in view of \eqref{5}, we obtain
\begin{equation}
\frac{{\rm d} R}{{\rm d}t}= \frac{1}{P_R}\,\frac{{\rm d} P}{{\rm d}t},
\label{6}
\end{equation}
where
\begin{equation}
P_R=\frac{\partial P(R,S)}{\partial R}=\frac{\alpha}{R}\left( \frac{RS}{S+1}\right)^{\alpha} + \frac{\gamma A}{R}\left( \frac{R}{S+1}\right)^{\gamma} >0.
\label{7}
\end{equation}

Taking into account \eqref{5} and \eqref{6}  as well as the divergence constraint \eqref{3}, we easily symmetrize the system of conservation laws \eqref{1} by rewriting it in the nonconservative form
\begin{equation}
\left\{
\begin{array}{l}
{\displaystyle\frac{1}{R P_R}\,\frac{{\rm d} P}{{\rm d}t} +{\rm div}\,{u} =0,\qquad
R\, \frac{{\rm d}u}{{\rm d}t}-({H},\nabla ){H}+{\nabla}q  =0 ,}\\[12pt]
{\displaystyle\frac{{\rm d}{H}}{{\rm d}t} - ({H} ,\nabla ){u} +
{H}\,{\rm div}\,{u}=0},\qquad
{\displaystyle\frac{{\rm d} S}{{\rm d} t} =0}.
\end{array}\right. \label{8}
\end{equation}
Equations \eqref{8} have absolutely the same form as the classical (one-fluid) MHD equations written in the nonconservative form. Writing down \eqref{8} in a matrix form, we get the symmetric system
\begin{equation}
A_0(U )\partial_tU+\sum_{j=1}^3A_j(U )\partial_jU=0
\label{9}
\end{equation}
which is {\it hyperbolic} ($A_0>0$) under assumptions \eqref{4}, where
\[
U= (P,u,H,S),\quad A_0= {\rm diag} (1/(RP_R) ,R ,R, R,1,1,1,1),
\]
\[
A_1=\left( \begin{array}{cccccccc} \frac{u_1}{RP_R}&1&0&0&0&0&0&0\\[3pt]
1&R u_1&0&0&0&{H_2}&{H_3}&0\\
0&0&R u_1&0&0&-{H_1}&0&0\\ 0&0&0&R u_1&0&0&-{H_1}&0\\
0&0&0&0&{u_1}&0&0&0\\
0&{H_2}&-{H_1}&0&0&{u_1}&0&0\\
0&{H_3}&0&-{H_1}&0&0&{u_1}&0\\ 0&0&0&0&0&0&0&u_1\\
\end{array} \right) ,
\]
\[
A_2=\left( \begin{array}{cccccccc} \frac{u_2}{RP_R}&0&1&0&0&0&0&0\\[3pt]
0&R u_2&0&0&-{H_2}&0&0&0\\ 1&0&R u_2&0&{H_1}&0&{H_3}&0\\ 0&0&0&R
u_2&0&0&-{H_2}&0\\
0&-{H_2}&{H_1}&0&{u_2}&0&0&0\\
0&0&0&0&0&{u_2}&0&0\\
0&0&{H_3}&-{H_2}&0&0&{u_2}&0\\0&0&0&0&0&0&0&u_2
\end{array} \right) ,
\]
\[
A_3=\left( \begin{array}{cccccccc} \frac{u_3}{RP_R}&0&0&1&0&0&0&0\\[3pt]
 0&R u_3&0&0&-{H_3}&0&0&0\\ 0&0&R u_3&0&0&-{H_3}&0&0\\ 1&0&0&R
u_3&{H_1}&{H_2}&0&0\\
0&-{H_3}&0&{H_1}&{u_3}&0&0&0\\
0&0&-{H_3}&{H_2}&0&{u_3}&0&0\\
0&0&0&0&0&0&{u_3}&0\\ 0&0&0&0&0&0&0&u_3
\end{array} \right)\; .
\]
The elementary symmetrization of conservations laws \eqref{1} giving the quasilinear symmetric hyperbolic system \eqref{8}/\eqref{9}, which looks like the symmetric system of the classical MHD equations, implies a number of results almost for nothing. For example, we have the local-in-time existence and uniqueness theorem \cite{Kato,VKh} in Sobolev spaces $H^s$ for the Cauchy problem for system \eqref{8}, with $s\geq 3$.

Following the idea of \cite{T05}, we can get a so-called {\it secondary symmetrization} of the symmetric hyperbolic system \eqref{9}. Such a secondary symmetrization was obtained in \cite{T05} by constructing a new energy integral for the magnetoacoustics system
\begin{equation}
A_0(\widehat{U} )\partial_tU+\sum_{j=1}^3A_j(\widehat{U} )\partial_jU=0
\label{10}
\end{equation}
associated with the nonlinear system \eqref{9}, where $\widehat{U}=(\widehat{P},\hat{u},\widehat{H},\widehat{S})$ is a constant vector. For the linear constant coefficients system \eqref{10} we have the standard conserved integral
\begin{equation}
\frac{{\rm d} I(t)}{{\rm d}t}=0
\label{11}
\end{equation}
in $\mathbb{R}^3$, with
\[
I(t)=\int\limits_{\mathbb{R}^3}(A_0(\widehat{U})U\cdot U )dx = \int\limits_{\mathbb{R}^3}\left(
\frac{1}{\widehat{R}\widehat{P}_R}P^2+\widehat{R}\,|u |^2 +|H |^2 +S^2\right){\rm d}x >0,
\]
where $\widehat{P}_R=P_R (\widehat{R},\widehat{S})$ and $\widehat{R}$ is a constant total density corresponding to the constant pressure $\widehat{P}$. We can write down the additional energy integral
\begin{equation}
\frac{{\rm d}J(t)}{{\rm d}t} =0,
\label{12}
\end{equation}
with
\[
J(t)=\int\limits_{\mathbb{R}^3}\left( \frac{1}{\widehat{P}_R}(\widehat{H}\cdot u )P -
\widehat{R}\,(u\cdot H )\right){\rm d}x,
\]
is a counterpart of the {\it cross-helicity} integral
\[
\frac{{\rm d}}{{\rm d}t}\int\limits_{\mathbb{R}^3}(u\cdot H )\,{\rm d}x =0
\]
in ideal {\it incompressible} MHD.

The combination of the energy integrals \eqref{11} and \eqref{12} gives the new conserved integral
\[
\frac{{\rm d}}{{\rm d}t}\bigl\{I(t)+2\lambda J(t)\bigr\} =\frac{{\rm d}}{{\rm d}t}\int\limits_{\mathbb{R}^3}({B}_0(\widehat{U})U\cdot U )\,{\rm d}x =0,
\]
where $\lambda$ is an arbitrary constant and the symmetric matrix
\[
{B}_0 =B_0(U)= \left( \begin{array}{cccccccc}
{\frac{1}{R P_R}} & {\frac{\lambda
H_1}{P_R}} & {\frac{\lambda H_2}{P_R}} &
{\frac{\lambda H_3}{P_R}} & 0 & 0 & 0 &0\\[9pt]
{\frac{\lambda H_1}{P_R}} & R & 0 & 0 &
-R\lambda & 0 & 0&0\\[6pt]
{\frac{\lambda H_2}{P_R}} & 0 & R & 0 & 0 &
-R\lambda & 0 &0\\[6pt]
{\frac{\lambda H_3}{P_R}}& 0 & 0 & R & 0 & 0 &
-R\lambda &0\\
0 &-R\lambda &0 & 0 & 1 & 0 & 0&0\\
0 & 0 &-R\lambda &0 & 0 & 1 & 0 &0\\
0 & 0 & 0 &-R\lambda &0 & 0 & 1&0 \\
0 & 0 & 0 &0 &0 & 0 & 0&1
\end{array}\right).
\]
We can easily find such a matrix ${\mathcal S}$ that ${\mathcal S}A_0=B_0$. This matrix is written down in \cite{T05,T09} (for our case, in \cite{T05,T09} one needs to replace $\rho$ and the square of the sound $c^2$ with $R$ and $P_R$ respectively). Let us now return to the nonlinear system \eqref{9} and assume that $\lambda =\lambda (U)$ be an arbitrary function of $U$.

Multiplying the symmetric hyperbolic system \eqref{9} from the left by the matrix ${\mathcal S}(U)$ and taking into account the divergence constraint \eqref{3}, we get the system
\begin{equation}
{\cal S}(U )A_0(U )\partial_t{U} +\sum_{j=1}^{3}{\cal S}(U )A_j(U )\partial_j{U} +
{\mathcal R}(U )\,{\rm div}\,{H} =0
\label{13}
\end{equation}
which is again symmetric for a suitable choice of the vector ${\mathcal R}(U)$ (we take ${\mathcal R}=-\lambda (1,0,0,0,H,0)$, see \cite{T05,T09}). Indeed, system \eqref{13} is written as
\begin{equation}
B_0(U )\partial_t{U} +\sum_{j=1}^{3}B_j(U )\partial_j{U}=0,
\label{14}
\end{equation}
with the symmetric matrices $B_j$ which are given in \cite{T05} (for our case we should again formally replace $\rho$ and $c^2$ appearing in the matrices $B_j$ in \cite{T05} with $R$ and $P_R$ respectively). System \eqref{14} is hyperbolic ($B_0>0$) under assumptions \eqref{4} supplemented by the following restriction on $\lambda$:
\begin{equation}
\lambda^2<\frac{P_R}{RP_R+|H|^2}.
\label{15}
\end{equation}

\section{Strong discontinuities in two-fluid MHD}
\label{s3}

We consider the two-fluid MHD equations \eqref{1} for $t\in [0,T]$ in the unbounded space domain $\mathbb{R}^3$ and suppose that $\Gamma (t)=\{ x_1-\varphi (t,x')=0\}$ is a smooth hypersurface in $[0,T]\times\mathbb{R}^3$, where
$x'=(x_2,x_3)$ are tangential coordinates. We assume that $\Gamma (t)$ is a surface of strong discontinuity for the conservation laws \eqref{1}, i.e., we are interested in solutions of \eqref{1} that are smooth on either side of $\Gamma (t)$. To be weak solutions of (\ref{1}) such piecewise smooth solutions should satisfy the Rankine-Hugoniot-type jump conditions
\begin{equation}
\left\{
\begin{array}{l}
[\mathfrak{j}_1]=0,\quad[\mathfrak{j}_2]=0,\quad [H_{N}]=0,\quad \mathfrak{j}\left[u_{N}\right] +
|N|^2[q]=0,\\[6pt]
\mathfrak{j}\left[{u}_{\tau}\right]=H_{N}[{H}_{\tau}],\quad
H_{N}[{u}_{\tau}]=\mathfrak{j}\left[{H}_{\tau}/R\right]
\end{array}
\right.
\label{16}
\end{equation}
at each point of $\Gamma$, where $[g]=g^+|_{\Gamma}-g^-|_{\Gamma}$ denotes the jump of $g$, with $g^{\pm}:=g$ in the domains
\[
\Omega^{\pm}(t)=\{\pm (x_1- \varphi (t,x'))>0\},
\]
and
\[
\mathfrak{j}_1^{\pm}=n (u_{N}^{\pm}-\partial_t\varphi),\quad
\mathfrak{j}_2^{\pm}=\rho (u_{N}^{\pm}-\partial_t\varphi),\quad \mathfrak{j}^{\pm}=R (u_{N}^{\pm}-\partial_t\varphi),\quad u_{N}^{\pm}=({u}^{\pm} {\cdot}{N}),
\]
\[
H_{N}^{\pm}=({H}^{\pm} {\cdot}{N}),\quad
{N}=(1,-\partial_2\varphi,-\partial_3\varphi ),\quad
{u}^{\pm}_{\tau}=(u^{\pm}_{\tau _1},u^{\pm}_{\tau _2}),\quad {H}^{\pm}_{\tau}= (H^{\pm}_{\tau _1}, H^{\pm}_{\tau _2}),
\]
\[
u^{\pm}_{\tau _i}=({u}^{\pm} {\cdot}{\tau}_i),\quad
H^{\pm}_{\tau
_i}=({H}^{\pm} {\cdot}{\tau}_i),\quad
{\tau}_1=(\partial_2\varphi,1,0),\quad
{\tau}_2=(\partial_3\varphi,0,1),\quad H_{N}|_{\Gamma}:=H_{N}^{\pm}|_{\Gamma};
\]
$\mathfrak{j}:=\mathfrak{j}^{\pm}|_{\Gamma}$ is the mass transfer flux across the discontinuity surface (the first two conditions in \eqref{16} imply $[\mathfrak{j}]=0$).

If we replace the first two conditions in \eqref{16} with $[\mathfrak{j}]=0$ and assume that $P=P(R)$, we obtain the Rankine-Hugoniot conditions in {\it isentropic} MHD with the fluid density $R$. That is, the classification of strong discontinuities in two-fluid MHD is closely related to that in isentropic MHD. However, there is a principal difference from isentropic MHD because in two-fluid MHD the equation of state $P=P(R,S)$ is two-parametric.  In particular, the continuity of the pressure on $\Gamma$ ($[P]=0$) does not imply $[R]=0$ and, similarly to full (non-isentropic) MHD, we can consider discontinuities with no flow across the discontinuity for which the pressure, the magnetic field and the velocity are continuous whereas the density (to be exact, the total density $R$ or the densities $n$ and $\rho$) may have a jump (see below).

From the mathematical point of view, there are two types of strong discontinuities: shock waves and characteristic discontinuities. Following Lax \cite{Lax57}, characteristic discontinuities, which are characteristic free boundaries, are called contact discontinuities. For the Euler equations of gas dynamics contact discontinuities are indeed contact  from the physical point of view, i.e., there is no flow across the discontinuity ($\mathfrak{j}=0$).

Recall that in classical MHD the situation with characteristic discontinuities is richer than in gas dynamics. Besides shock waves there are three types of characteristic discontinuities \cite{BThand,LL}: current-vortex sheets, contact discontinuities and Alfv\'{e}n discontinuities. Current-vortex sheets and MHD contact discontinuities are contact from the physical point of view, but Alfv\'{e}n discontinuities are not. In fact, the same situation we have in two-fluid MHD. As in one-fluid non-isentropic MHD, we consider the following four types of discontinuities:
\begin{itemize}
  \item[1)] shock wave ($\mathfrak{j}\neq 0$ and $[R ]\neq 0$),
  \item[2)] current-vortex sheet ($\mathfrak{j}=0$ and $H_{N}|_{\Gamma}=0$),
  \item[3)] contact discontinuity ($\mathfrak{j}=0$ and $H_{N}|_{\Gamma}\neq 0$),
  \item[4)] Alfv\'{e}n discontinuity ($\mathfrak{j}\neq 0$ and $[R ]= 0$).
\end{itemize}
Let us first consider shock waves.

\paragraph{Shock waves} Since $\mathfrak{j}\neq 0$, we have $u_{N}^{\pm}|_{\Gamma}- \partial_t\varphi\ne 0$. Then, the first two conditions in \eqref{16} imply
\begin{equation}
[S]=0.
\label{17}
\end{equation}
The rest jump conditions are the same as in isentropic MHD. The continuity of the ``entropy'' \eqref{17} means that shock waves for the conservation laws \eqref{1} are indeed similar to shock waves in isentropic MHD (see the next section for more details).

\paragraph{Current-vortex sheets} Assuming that $\mathfrak{j}=0$ and $H_{N}|_{\Gamma}=0$, from \eqref{16} we get the boundary conditions
\begin{equation}
[q]=0,\quad H_{N}^{\pm}=0,\quad \partial_t\varphi =u_N^{\pm}\quad \mbox{on}\ \Gamma (t).
\label{18}
\end{equation}
At the same time, the tangential components of the velocity and the magnetic field may undergo any jump: $[u_{\tau}]\neq 0$ and $[H_{\tau}]\neq 0$. Clearly, the same is true for the densities $n$ and $\rho$. In principle, one of them may be continuous but, in general,
$[n]\neq 0$ and $[\rho ]\neq 0$.

\paragraph{Contact discontinuities} For a contact discontinuity, it follow from \eqref{16} that
\begin{equation}
[P]=0,\quad [u]=0,\quad [H]=0,\quad \partial_t\varphi =u_N^{\pm}\quad \mbox{on}\ \Gamma (t).
\label{19}
\end{equation}
Moreover, the densities $n$ and $\rho$ (and the total density $R$) may undergo any jump.

\paragraph{Alfv\'{e}n discontinuities}
It follows from $[\mathfrak{j}]=0$ and $[R]=0$ that $[u_N]=0$. Then, the first two conditions in \eqref{16} imply $[n]=[\rho]=0$ and, hence, $[P]=0$. In view of $[R]=0$, the last two (vector) conditions in \eqref{16} form the linear homogeneous system for the jumps $[u_{\tau}]$ and $[H_{\tau}]$. If $\mathfrak{j}\neq\pm H_N\sqrt{R}$ at each point of $\Gamma$ (here $R:=R^{\pm}|_{\Gamma}$), then $[u_{\tau}]\equiv 0$ and $[H_{\tau}]\equiv 0$. In this case we have no discontinuity at all: $[P]=[n]=[\rho]=0$ and $[u]=[H]=0$. As in classical MHD, we assume that $\mathfrak{j} = H_N\sqrt{R}$ at each point of $\Gamma$. Taking this into account, we finally obtain the boundary conditions
\begin{equation}
[P]=0,\quad [S]=0,\quad [H_N]=0,\quad [|H|^2]=0,\quad j=H_N\sqrt{R},\quad [u]=\frac{[H]}{\sqrt{R}} \qquad \mbox{on}\ \Gamma (t)
\label{20}
\end{equation}
for the components of the vectors $U^{\pm}$ ($U^{\pm}:=U$ in $\Omega^{\pm}(t)$). The boundary conditions \eqref{20} on a surface of Alfv\'{e}n discontinuity in two-fluid MHD are absolutely the same as that for Alfv\'{e}n discontinuities in classical MHD \cite{BThand,LL} describing the motion of a fluid having the density $R(t,x)$ and the entropy $S(t,x)$.

\section{Structural stability conditions}
\label{s4}

Basing on the known results for classical MHD, we now discuss structural stability conditions for initial data for the free boundary problems for all the types of strong discontinuities in two-pase MHD introduced in Section \ref{s3}. These conditions should guarantee the local-in-time existence in Sobolev spaces of solutions of these problems (usually they also provide the uniqueness of a solution). We begin with shock waves.

\subsection{Shock waves}

The free boundary problem for shock waves is the problem for the systems
\begin{equation}
A_0(U^{\pm})\partial_tU^{\pm}+\sum_{j=1}^3A_j(U^{\pm} )\partial_jU^{\pm}=0\quad \mbox{in}\ \Omega^{\pm}(t)
\label{21}
\end{equation}
(cf. \eqref{9}) with the boundary conditions
\begin{equation}
\left\{
\begin{array}{l}
[\mathfrak{j}]=0,\quad [S]=0,\quad [H_{N}]=0,\quad \mathfrak{j}\left[u_{N}\right] +
|N|^2[q]=0,  \\[6pt]
\mathfrak{j}\left[{u}_{\tau}\right]=H_{N}[{H}_{\tau}],\quad
H_{N}[{u}_{\tau}]=\mathfrak{j}\left[{H}_{\tau}/R\right]
\end{array}
\right.
\label{22}
\end{equation}
on $\Gamma (t)$ (cf. \eqref{16}, \eqref{17}) and corresponding initial data for $U^{\pm}$ and $\varphi$ at $t=0$. We can reduce this problem to that in the fixed domains $\mathbb{R}^3_{\pm}=\{\pm x_1>0,\quad x'\in\mathbb{R}^2\}$ by the simple change of variables
\begin{equation}
\tilde{x}_1=x_1-\varphi (t,x'),\quad \tilde{x}'=x'.
\label{chv}
\end{equation}
Dropping tildes, from systems \eqref{21} we obtain
\begin{equation}
A_0(U^{\pm})\partial_tU^{\pm}+\widetilde{A}_1(U^{\pm},\varphi )\partial_1U^{\pm}+A_2(U^{\pm} )\partial_2U^{\pm}+A_3(U^{\pm} )\partial_3U^{\pm}=0\quad \mbox{for}\ x\in \mathbb{R}^3_{\pm},
\label{23}
\end{equation}
where $\widetilde{A}_1$ is the so-called boundary matrix:
\[
\widetilde{A}_1=\widetilde{A}_1(U,\varphi )= A_1(U)-A_0(U)\partial_t\varphi -A_2(U)\partial_2\varphi -A_3(U)\partial_3\varphi .
\]
The boundary conditions for \eqref{23} are \eqref{22} on the plane $x_1=0$.

It is well-known that the necessary condition for the well-posedness of the above problem is that the number of boundary conditions should be one unit greater than the number of incoming characteristics of the 1D counterparts (with $A_2=A_3=0$) of systems \eqref{23} for fixed (``frozen'') $U^{\pm}$ and $\varphi$ satisfying the boundary conditions (roughly speaking, one of the boundary conditions is needed for finding the unknown function $\varphi$). Since the number of incoming characteristics is defined by the number of positive/negative eigenvalues of the matrices $\left(A_0(U^{\pm})\right)^{-1}\widetilde{A}_1(U^{\pm},\varphi )$, for shock waves (for them these matrices have no zero eigenvalues) this is equivalent to the Lax's $k$-shock conditions
\[
\lambda_{k-1}^- <\partial_t\varphi <\lambda_k^-,\quad \lambda_k^+ <\partial_t\varphi <\lambda_{k+1}^+
\]
for some integer $k$, where for our case of system of eight equations $1\leq k\leq 8$ and $\lambda_j^\pm$ ($j=\overline{1,8}$) are the eigenvalues of the matrices
\[
A_N^{\pm}:=\left(A_0(U^{\pm})\right)^{-1}\left( A_1(U^{\pm})-A_2(U^{\pm})\partial_2\varphi -A_3(U^{\pm})\partial_3\varphi \right)|_{x_1=0},
\]
with some fixed $U^{\pm}$ and $\varphi$ satisfying the boundary conditions \eqref{22} at $x_1=0$. Moreover, $\lambda_j^\pm$ are numbered as
\[
\lambda_1^-\leq \ldots \leq\lambda_8^-,\quad \lambda_1^+\leq \ldots \leq\lambda_8^+,
\]
and we take $\lambda^-_0:=\partial_t\varphi/2$ and $\lambda_9^+:=2\partial_t\varphi$.

Since system \eqref{8} have absolutely the same form as the classical MHD equations, we already know the eigenvalues $\lambda_j^\pm$ (see, e.g., \cite{BThand,LL}):
\begin{equation}
\begin{split}
& \lambda_1^{\pm}=u_N^{\pm}-c_f^{\pm},\quad \lambda_2^{\pm}=u_N^{\pm}-c_a^{\pm},\quad \lambda_3^{\pm}=u_N^{\pm}-c_s^{\pm}, \quad
\lambda_4^{\pm}=\lambda_5^{\pm}=u_N^{\pm},\\
& \lambda_6^{\pm}=u_N^{\pm}+c_s^{\pm},\quad \lambda_7^{\pm}=u_N^{\pm}+c_a^{\pm},\quad \lambda_8^{\pm}=u_N^{\pm}+c_f^{\pm},
\end{split}
\label{24}
\end{equation}
where $c_a^{\pm} = H_N^{\pm}/\sqrt{R^{\pm}}$ are the Alfv\'{e}n speeds (ahead and behind of the shock) in the normal direction,
\[
c_s^{\pm}=\frac{1}{\sqrt{2}}\sqrt{(c^{\pm})^2+(c_A^{\pm})^2-\sqrt{\left((c^{\pm})^2+(c_A^{\pm})^2\right)^2-4(c^{\pm}c_a^{\pm})^2}}
\]
are the slow magnetosonic speeds,
\[
c_f^{\pm}=\frac{1}{\sqrt{2}}\sqrt{(c^{\pm})^2+(c_A^{\pm})^2+\sqrt{\left((c^{\pm})^2+(c_A^{\pm})^2\right)^2-4(c^{\pm}c_a^{\pm})^2}}
\]
are the fast magnetosonic speeds,  $c^{\pm}=\left(P_R(R^{\pm},S^{\pm})\right)^{1/2}$ play the role of sound speeds (ahead and behind of the shock) in two-fluid MHD, and $c_A^{\pm} = |H^{\pm}|/\sqrt{R^{\pm}}$. Note also that all the values in \eqref{24} are written at $x_1=0$.
Clearly, as in classical MHD, we have two types of admissible shocks satisfying the Lax's $k$-shock conditions: {\it fast shock waves} and {\it slow shock waves}. Without loss of generality we assume that $u_N^{+}|_{x_1=0}>\partial_t\varphi$. Then, fast shock waves are 1-shocks satisfying the inequalities
\[
c_f^-<u_N^--\partial_t\varphi,\quad c_a^+<u_N^+-\partial_t\varphi <c_f^+,
\]
and slow shock waves are 3-shocks satisfying the inequalities
\[
c_s^-<u_N^--\partial_t\varphi <c_a^-,\quad u_N^+-\partial_t\varphi <c_s^+.
\]
We note that fast shock waves are {\it extreme shocks} in the sense that ahead of the shock there are no incoming waves.

According to the results in \cite{M1,M2,Met} and their extension to hyperbolic symmetrizable systems with characteristics of variable multiplicities \cite{Kwon,MZ}, all uniformly stable shocks are structurally stable. Roughly speaking (we do not discuss regularity, compatibility conditions, etc.), this means that if the uniform Lopatinski condition holds at each point of the initial shock, then this shock exists locally in time. In other words, as soon as planar shock waves are uniformly stable according to the linear analysis with constant coefficients, we can make the conclusion about structural stability of corresponding nonplanar shocks. In this sense, the linear analysis with constant coefficients plays crucial role for shock waves.

Without loss of generality we can consider the unperturbed shock wave with the equation $x_1=0$. Then, the linearization of \eqref{22}, \eqref{23} about a constant solution, with $\widehat{U}^{\pm}=(\widehat{P}^{\pm},\hat{u}^{\pm},\widehat{H}^{\pm},\widehat{S}^{\pm})={\rm const}$ and $\hat{\varphi} =0$, gives a linear constant coefficients problem. For the perturbation of the ``entropy'' (we again denote it by $S$) we obtain the separate problem (cf. \eqref{8}, \eqref{17})
\[
\partial_tS^{\pm} +(\hat{u}^{\pm}\cdot S) =f^{\pm}\quad \mbox{for}\ x\in\mathbb{R}_3^{\pm},\qquad [S]=g \quad \mbox{at}\ x_1=0,
\]
where, as usual, we introduce artificial source terms (for the above equations they are $f^{\pm}(t,x)$ and $g(t,x')$) to make the linear problem inhomogeneous. For this problem we easily deduce the energy identity
\begin{equation}
I(t)-\int\limits_0^t\int\limits_{\mathbb{R}^2}
\left.\left( [\hat{u_1}](S^-)^2 +2\hat{u}^+_1S^-g+\hat{u}^+_1g^2 \right)\right|_{x_1=0}\,{\rm d}x'{\rm d}s=I(0)+2\sum_{\pm}\int\limits_0^t\int\limits_{\mathbb{R}^3_{\pm}}f^{\pm}S^{\pm}{\rm d}x{\rm d}s,
\label{27}
\end{equation}
where
\[
[\hat{u}_1]=\hat{u}_1^+-\hat{u}_1^-\quad\mbox{and}\quad
I(t)=\sum_{\pm}\|{S}^{\pm}(t)\|^2_{L^2(\mathbb{R}^3_{\pm})}.
\]
For {\it compressive} shocks ($[R]>0$), it follows from the first condition in \eqref{22} that $[u_N]<0$. That is, for the unperturbed planar compressive shock we have $[\hat{u}_1]<0$. Then, for compressive shocks, by standard simple arguments, from \eqref{27} we get the a priori estimate without loss of derivatives
\[
\begin{split}
\sum_{\pm}\Big(\|{S}^{\pm}\|_{L^2([0,T]\times\mathbb{R}^3_{\pm})}+ & \|S^{\pm}_{|x_1=0}\|_{L^2([0,T]\times\mathbb{R}^2)}\Big) \\  & \leq C\left\{ \sum_{\pm}\left(\|{S}^{\pm}_{|t=0}\|_{L^2(\mathbb{R}^3_{\pm})}
+\|f^{\pm}\|_{L^2([0,T]\times\mathbb{R}^3_{\pm})}\right) +\|g\|_{L^2([0,T]\times\mathbb{R}^2)}\right\},
\end{split}
\]
where $C>0$ is a constant.

Taking into account the above a priori estimate, the linearized problem satisfies the uniform Lopatinski condition as soon as its subproblem for the perturbations of $P^{\pm}$, $u^{\pm}$ and $H^{\pm}$ does. This subproblem totally coincides with the linearized constant coefficients problem for shock waves in isentropic MHD with the equation of state
\[
{P}({R})=c_1{R}^{\alpha}+ c_2{R}^{\gamma},
\]
for the unperturbed constant solution, where $c_1=({S}/({S}+1))^{\alpha}\geq 0$ and $c_2=A(1/({S}+1))^{\gamma}>0$ are constants for a fixed $S=\widehat{S}^{+}=\widehat{S}^-$ (recall that $[\widehat{S}]=0$). For $\alpha =\gamma$ this equation of state is that for a polytropic gas (if we are speaking about one-fluid isentropic MHD with the fluid density $R$). But, even if $\alpha \neq\gamma$, in the above equation of state $P(R)$ is a convex function (for $\alpha\geq 1$ and $\gamma >1$).  Without magnetic field it was shown by Majda \cite{M1} that all isentropic {\it compressive} shock waves are uniformly stable provided that the equation of state is convex. Their local-in-time existence was proved in \cite{M2} (see also \cite{Met} for some improvements of the results in \cite{M2}). In the limit of weak magnetic field ($H\rightarrow 0$), the uniform stability of all compressive {\it fast} shock waves in isentropic MHD (implying their local-in-time existence) was proved by M\'etivier and Zumbrun \cite{MZ} for equations of state satisfying the uniform Lopatinski conditon in isentropic gas dynamics \cite{M1}, in particular, for convex equations of state. Basing on these results, we can thus make the conclusion that all compressive fast shock waves in two-fluid MHD are structurally stable.

Note that in full (non-isentropic) MHD the uniform stability of fast shock waves was proved by Blokhin and Trakhinin (see \cite{BThand} and references therein) by the energy method provided that the equation of state satisfies the uniform Lopatinski conditon in (full) gas dynamics.
Regarding the general (not weak) strength of magnetic field, a complete 2D stability analysis of fast MHD shock waves was carried out in \cite{T} for the polytropic gas equation of state. Taking into account the results of \cite{MZ,Kwon} for hyperbolic symmetrizable systems with characteristics of variable multiplicities (this class contains the MHD system), uniformly stable fast shock waves found in \cite{T} exist locally in time. For slow shock waves, some results about their stability can be found in \cite{BThand,Fil}. Since in isentropic MHD we do not have analogues of the results in \cite{BThand,Fil,T} obtained for shock waves in full MHD and since such results for isentropic MHD shocks could be directly carried over two-fluid MHD, there is a big motivation of the stability analysis for shock waves in isentropic MHD under an arbitrary strength of magnetic field.

\subsection{Characteristic discontinuities}

For {\it current-vortex sheets} in ideal compressible MHD the crucial idea was the secondary symmetrization of the MHD system proposed in \cite{T05}. The usage of this symmetrization gives a condition sufficient for neutral stability of a planar current-vortex sheet, i.e. a condition sufficient for the fulfillment of the Lopatinski condition for the constant coefficients linearized problem. Having in hand this condition and assuming that it holds at each point of the initial current-vortex sheet, the local-in-time existence and uniqueness theorem for the nonlinear problem was independently proved in \cite{ChW1,ChW2} and \cite{T09}.

As we have seen in Section \ref{s3}, we have the same secondary symmetrization  for two-fluid MHD (see \eqref{13}, \eqref{14}). We can just reformulate the sufficient stability condition from \cite{T05,T09} for current-vortex sheets in two-pase MHD:
\begin{equation}
G(U^+,U^-)|_{\Gamma}>0,
\label{28}
\end{equation}
where
\[
G(U^+,U^-)=|\sin(\psi^+-\psi^-)|
\min\left\{\frac{\beta^+}{|\sin{\psi^-}|}\,,\,
\frac{\beta^-}{|\sin{\psi^+}|}\right\}-|[u']|,
\]
\[
u'^{\pm}=(u_2^{\pm},u_3^{\pm}),\quad H'^{\pm}=(H_2^{\pm},H_3^{\pm}),\quad
\beta^{\pm}=\frac{{c}^{\pm}{c}^{\pm}_{A}}{\sqrt{({c}^{\pm})^2+({c}^{\pm}_{A})^2}}\,,\quad
\cos\psi^{\pm}=\frac{([u']\cdot H'^{\pm})}{ |[u']|\,|H'^{\pm}|},
\]
and the speeds $c^{\pm}$ and $c_A^{\pm}$ were defined just after \eqref{24}. Here it is assumed that the vectors $H'^+$ and $H'^-$ are not collinear at each point of $\Gamma$. We now make the conclusion that current-vortex sheets in two-fluid MHD exist locally in time provided that condition \eqref{28} holds for $t=0$.  Of course, we should also assume that the initial data satisfy the hyperbolicity conditions \eqref{4}, the non-collinearity condition for the magnetic fields $H'^{\pm}$ and appropriate compatibility conditions.

Let us now consider  {\it contact discontinuities} (see \eqref{19}). In classical MHD with the fluid density $\rho$, etc., the important assumption for them made in \cite{MTT1,MTT2} was the continuity of the product of the density and the square of the sound speed: $[\rho c^2]=0$. In particular, this is true for a polytropic gas. The important consequence of this assumption is that $[\partial_1u]=0$ for solutions of the free boundary problem reduced to that in the fixed domains $\mathbb{R}^{\pm}_3$. For our case of two-fluid MHD, this assumption reads:
\[
[RP_R]=\left[ \alpha\left( \frac{RS}{S+1}\right)^{\alpha} + \gamma A\left( \frac{R}{S+1}\right)^{\gamma}\,\right]=0
\]
(see \eqref{7}). If $\alpha=\gamma$, then $[RP_R]=\gamma [P]=0$ (cf. \eqref{19}). Since \eqref{8} has absolutely the same form as the classical MHD system, as in  \cite{MTT1}, we can prove that $[\partial_1u]=0$. The rest arguments towards the proof of a counterpart of the local-in-time existence and uniqueness theorem for contact discontinuities are totally the same as in  \cite{MTT1,MTT2}. We only note that such a theorem was proved in  \cite{MTT2} for the 2D case provided that the Rayleigh-Taylor sign condition $[\partial P/\partial N] < 0$ on the jump of the normal derivative of the pressure is satisfied at each point of the initial discontinuity.

We can thus make the conclusion that for the equation of state \eqref{2}, with $\alpha =\gamma$, 2D contact discontinuities in two-fluid MHD exist locally in time under the Rayleigh-Taylor sign condition $[\partial P/\partial N] < 0$ satisfied for $t=0$.

At last, we briefly discuss {\it Alfv\'{e}n discontinuities}. In classical MHD, there are still no nonlinear (structural stability) results for them. At the same time, the domains of neutral stability and violent instability of planar Alfv\'{e}n discontinuities were numerically found in \cite{IT}. It is still unclear whether neutrally stable nonplanar Alfv\'{e}n discontinuities do exist locally in time. However, the linear results in \cite{IT} are automatically carried over two-fluid MHD because the free boundary problem \eqref{20}, \eqref{21} has absolutely the same form as that for Alfv\'{e}n discontinuities in classical MHD.

\bigskip

\section*{Acknowledgements}
The research of L.Z. Ruan was supported in part by the Natural Science Foundation of China $\#$11771169, $\#$11331005, $\#$11301205, $\#$11871236, Program for Changjiang Scholars and Innovative Research Team in University $\#$IRT17R46, and the Special Fund for Basic Scientific Research of Central Colleges $\#$CCNU18CXTD04.

\end{document}